\documentclass[12pt]{amsart}
\usepackage{amssymb,latexsym,amsmath,amscd,mathrsfs,yfonts}
\usepackage{fullpage}
\input xy
\xyoption{all}

\newcommand{\ZZ}{\mathbb Z}
\newcommand{\QQ}{\mathbb Q}
\newcommand{\CC}{\mathbb C}

\newcommand{\G}{\mathcal G}

\newcommand{\LL}{\mathbb L}
\newcommand{\Aff}{\mathbb A}

\def\cA{\mathcal{A}}
\def\cB{\mathcal{B}}

\def\QCoh{\textup{QCoh}}

\def\per{\textup{per}}
\def\Perf{\mathtt{Perf}}
\def\Sch{\mathtt{Sch}}

\def\dim{\textup{dim}}
\def\Spec{\textup{Spec}}

\def\innHom{\Sheaf{H}om}

\def\Hom{\textup{Hom}}

\def\ker{\textup{ker}}

\def\op{\textup{op}}
\def\pt{\textup{pt}}

\def\rep{\textup{rep}}

\def\Hqe{\mathtt{Hqe}}

\def\Var{\mathtt{Var}}
\def\Top{\mathtt{Top}}
\def\NCS{\mathtt{NCS}}
\def\NCCorr{\mathtt{NCC}}
\def\NCM{\mathtt{NCM}}
\def\DG{\mathtt{DGcat}}
\def\MC{\mathtt{NCM_{\text{CY}}}}
\def\H{\textup{H}}
\def\Corr{\mathtt{Corr}}
\def\Obj{\textup{Obj}}
\def\DM{\mathcal{D}_{\mathcal{M}}}
\def\D{\mathscr D}
\def\cC{\mathcal C}

\def\cD{\mathcal D}

\def\TT{\mathbb{T}}

\def\to{\rightarrow}
\newcommand{\map}{\longrightarrow}

\newcommand{\sheaf}{\mathcal{O}}

\newcommand{\Sheaf}{\mathscr}
\newcommand{\beq}{\begin{eqnarray}}
\newcommand{\beqn}{\begin{eqnarray*}}
\newcommand{\eeq}{\end{eqnarray}}
\newcommand{\eeqn}{\end{eqnarray*}}

\newtheorem{thm}{Theorem}[section]

\newtheorem{prop}[thm]{Proposition}

\newtheorem{ex}{Example}
\newtheorem{defn}[thm]{Definition}
\newtheorem{rem}[thm]{Remark}

\begin{document}

\title{Noncommutative geometry in the framework of differential graded categories}
\author{Snigdhayan Mahanta}

\maketitle

\setcounter{tocdepth}{2}
\tableofcontents


\section{Introduction} \label{DerivedNCG}

This expository note is meant to be a brief introduction to noncommutative geometry in a differential graded (DG) framework, {\it i.e.,} such categories playing the role of spaces. Keeping in mind the significance of certain correspondences in classical geometry, we are naturally led to include them as some sort of generalised morphisms of spaces. In this sense this geometry is {\it motivic}. However, although the construction of the category of noncommutative spaces will follow closely that of motives, the resulting category of noncommutative spaces will not even be additive (only semiadditive). Passing on to the homotopy categories of the DG categories (considered as spaces) one recovers most of the  results known at the level of triangulated categories. However, in this setting one does not run into some unpleasant technical problems which one would otherwise have to deal with at the level of triangulated categories. The theory seems to be a blend between algebraic topology (or homotopical geometry) and algebraic geometry. There is a possibility of recasting many different models of noncommutative geometry in this general setting. We also include some pointers to some other areas of mathematics, which are well adapted to be seen in this context. 

It must also be emphasized that this is noncommutative geometry 'at a large scale' (after Ginzburg \cite{GinNCG}) and, therefore, there are some natural new phenomena which are not quite compatible with the classical picture. One such instance is the isomorphism between certain classical spaces, {\it e.g.,} an abelian variety and its dual, which need not be isomorphic as classical spaces (varieties). For noncommutative geometry 'at a small scale', {\it i.e.,} viewed as a deformation of classical geometry one may look at, {\it e.g.,} \cite{Kap},\cite{Lau3}.

The outline of the construction presented here can be found in a recent preprint of
  Kontsevich \cite{KonMot}. In fact, this is a simplified version of the proposed one in {\it ibid.} Readers should also refer to the articles of J. Lurie and To\"en-Vezzosi \cite{DAG2,DAG3,TVDAG1,TVDAG2}, which seem to have developed a geometry based on a functor of points  approach from simplicial algebras (equivalent to connective differential graded algebras) to simplicial sets. The form of geometry in their parlance is homotopical algebraic geometry (HAG) or derived algebraic geometry (DAG)\footnote{The author would like to thank D. Ben-Zvi for providing quite convincing arguments to dispel any idea that HAG or DAG can subsume noncommutative geometry}. The material presented here is mostly modelled on the ICM talk of Keller \cite{KelDG}, though there are some minor deviations.

\vspace{2mm}
\noindent
{\bf Acknowledgements:} The author would like to thank Matilde Marcolli and Yu. I. Manin for generously sharing their thoughts and providing useful comments. A special word of thanks goes out to B. Keller for suggesting several changes and corrections based on an initial draft. The author would also like to thank  B. To\"en for patiently answering many questions and the anonymous referee for numerous meticulous remarks. The author is grateful to the Max-Planck-Institut f\"ur Mathematik, Bonn and the Fields Institute, Toronto for the hospitality, while this work was carried out.

\section{Noncommutative geometry in a DG framework}

For a long time it was felt that the language of triangulated
categories is deficient for many purposes in geometry. The language of DG
categories seems to have resolved most of the technical and aesthetic
problems. We first prepare the readers for the seemingly abstruse
definition of the category of noncommutative spaces. We propose a theory over an algebraically closed field of characteristic zero and by {\it Lefschetz Principle} there is no harm in assuming our ground field $k$ to be actually $\CC$. This reduces a lot of technical difficulties. For brevity, we denote most of the functors by their underived notation, for instance, ${\otimes}^{\mathbb{L}}$ is written simply as $\otimes$.

\subsection{Motivation} 

The traditional way of doing geometry with the emphasis on spaces is deficient
  in many physical situations. Most notably, due to Heisenberg's Uncertainty
  Principle one is forced to consider polynomial algebras with noncommuting
variables, like Weyl algebras. One has to do away with the notion of points of
a space quite naturally. However, one has perfectly well-defined algebras,
albeit noncommutative, with which one can work. One very successful approach
from this point of view is that of Connes \cite{ConBook}. It has many
applications and a large part of the classical (differential) geometry can be subsumed in this setting. One might also want to take a
closer look at the key features of classical (algebraic) geometry and try to
generalise them.

\bigskip
\noindent
{\bf From spaces to categories; from functions to sheaves:}
It is quite common in mathematics to study an object via its representations (in
an appropriate sense). It is neat to assemble all representations into a
category and study it. In this manner from groups one is led to the study of 
Tannakian categories, from algebras to that of certain triangulated categories and so on. This process is roughly some sort of categorification. 

We have already done away with the traditional notion of a space and its points. For the time being it is described by its functions. The topology of a space allows
us to define functions locally and glue them (if possible to a global one). A better way of keeping track of such information is using the language of the sheaf of local sections or functions on a space. Every classical space comes hand in hand with its structure sheaf of admissible functions {\it e.g.,}
continuous, smooth, holomorphic, algebraic, etc. according to the structure of
the underlying space. The representations of the structure sheaf, which for us are
nothing but quasicoherent sheaves, determine the space. In this manner one
replaces the notion of a space by its category of quasicoherent sheaves, an
idea that goes back to Gabriel, Grothendieck, Manin and Serre. 

The category of quasicoherent sheaves is a Grothendieck category when the scheme is quasicompact and quasiseparated \cite{ThoTro}. There are many approaches towards developing a theory by treating abelian categories (or some modifications
thereof, like Grothendieck categories) as the category of quasicoherent
sheaves on noncommutative spaces, {\it e.g.,} \cite{AZ,VDB1,Ros2}.

\begin{rem}
There is another point of view inspired by the Geometric Langlands programme
and the details can be found, for instance, in \cite{Fre}. The guiding principle
here is a generalisation of
Grothendieck's {\it faisceaux-fonctions correspondence}. The faisceaux-fonctions correspondence
appears naturally in the context of \'etale $\ell$-adic
sheaves. Associated to any complex of \'etale $\ell$-adic sheaves
$\Sheaf{K}^\bullet$ over a variety $V$ defined over a finite field
$\mathbb{F}_q$ is a function $f^{\Sheaf{K}^\bullet} : V\map \CC$ given
by 

\beqn
f^{\Sheaf{K}^\bullet}(x) = \sum (-1)^{i} \textup{Tr} (\textup{Fr}_{\bar{x}}\,|\,\textup{H}^i
(\Sheaf{K}^\bullet)_{\bar{x}}). 
\eeqn Here $x\in V(\mathbb{F}_q)$ and $\bar{x}$ denotes a geometric point of
$V$ over $x$. Of course, one has to fix an identification $\overline{\QQ}_\ell\overset{\sim}{\map} \CC$. According to
Grothendieck all {\it interesting} functions appear in this manner and
extrapolating this idea we regard constructible sheaves as the only source of
interesting functions over $\CC$. 

The lack of Verdier Duality, which is a generalisation of Poincar\'e Duality
and hence an important feature, makes the na\"ive category of constructible
$\ell$-adic sheaves undesirable. Instead one works with the category of
so-called perverse sheaves.\footnote{It is known that the derived category of
  coherent sheaves also admits a dualising complex imitating Grothendieck--Serre duality in place of Verdier duality (see Proposition 1 \cite{Bez}).} They are objects which live in a bigger derived category. 
Via a version of the Riemann-Hilbert correspondence over $\CC$ the category of
perverse sheaves (of middle perversity) is
equivalent to the category of regular holonomic $\D$-modules. More precisely,
let $X$ be a complex manifold, $D^b_{rh}(\D_X)$ denote the bounded derived category
of complexes of $\D_X$-modules with regular holonomic cohomologies and
$D^b_c(\CC_X)$ denote the bounded derived category of sheaves of complex
vector spaces with constructible cohomologies. Then Kashiwara proved in
\cite{Kas84} $\Sheaf{RH}om_{\D_X}(-,\sheaf_X):D^b_{rh}(\D_X)\overset{\sim}{\to}
D^b_c(\CC_X)^{\op}$ is an equivalence of triangulated categories. Under this
equivalence the standard t-structure on $D^b_{rh}(\D_X)$, whose heart is the abelian category of
regular holonomic $\D$-modules on $X$, is mapped to the heart of the
t-structure of middle perversity on $D^b_c(\CC_X)$. The heart of this
t-structure is the category of perverse sheaves (of middle perversity), which can be regarded as
another generalisation of functions. As opposed to a quasicoherent sheaf, the
model for a function in this setting is a $\D$-module, which is roughly a quasicoherent sheaf with a flat connection. A quasicoherent sheaf (resp. a $\D$-module) corresponds to a polynomial (resp. a constructible locally constant) function.
\end{rem}

\bigskip
\noindent
{\bf The passage to derived categories:} In the category of smooth schemes any morphism $f:X\to Y$
gives rise to two canonical functors on the category of sheaves, {\it viz.},
pull-back $f^\ast$ and push-forward $f_\ast$. One should naturally expect any
generalisation of classical geometry to allow such operations. We see that
restricting to abelian categories is not enough as functors like push-forwards
are not exact. The natural framework for such functors to exist is that of
derived categories or abstract triangulated categories. Besides, if one
chooses to work with perverse sheaves as substitutes for functions one has to view them as elements  of an abelian category sitting inside a bigger derived category. 

\bigskip
\noindent
{\bf Adding correspondences to morphisms:} Denoting by $\Var$ the
category of complex algebraic varieties, $\Top$ that of {\it nice} topological spaces (here {\it nice} should imply all
properties typical of the complex points of a complex algebraic variety) one 
has a tensor functor $\Var \to \Top$ associating to a complex
algebraic variety its underlying space with analytic topology. The tensor
structure on the two categories is given by direct product. To a topological
space in $\Top$ one can associate its singular cochain complex which is also a
tensor functor to $D_{ab}$, the category of complexes of finitely generated abelian groups whose  cohomology is bounded. According
to Be{\u\i}linson and Vologodsky \cite{BeiVol} the basic objective of the theory
of motives is to fill in a commutative diagram
\beqn
\begin{CD}
\Var @>>> \DM \\
@VVV      @VVV \\
\Top @>>> D_{ab}
\end{CD}
\eeqn where $\DM$ is the rigid tensor triangulated category of motives. The
upper horizontal arrow should be faithful and defined purely
geometrically and the right vertical arrow should respect the tensor
structures. In order to construct the upper horizontal arrow one first needs to enrich $\Var$ to include correspondences (modulo some equivalence relation). This endows $\Var$ with an
additive structure. 

\bigskip
\noindent
{\bf Triangulated structure is not enough:} The goal is to
construct a rigid tensor category of {\it motivic} noncommutative
spaces which allows basic operations like pull-back, push-forward and finite
correspondences (as morphisms). In the classical setting, we have a
construction of $\DM$ as a triangulated category due to Voevodsky (see {\it e.g.,} \cite{FSV}). However,
one would like to extract the {\it right} category of motives inside it
(possibly as an abelian rigid tensor category). One basic operation is direct product, which endows $\Var$ with the tensor structure. It should also survive in $\DM$. The tensor product of two triangulated categories unfortunately does not carry a natural triangulated structure. Also one
runs into trouble in trying to define inner $\Hom$'s. This is where the framework of DG
(differential graded) categories comes in handy.

\subsection{Overview of DG categories} \label{DGdiscussion}

Before we are able to spell out the definition of the category of
noncommutative spaces we need some preparation on DG categories, which will be
quite concise. For details we refer the readers to {\it e.g.,}
\cite{Dri},\cite{KelDG},\cite{dgToe}. They can be defined  over $k$, where $k$ is not necessarily a
field. However, as mentioned before, we set $k=\CC$  and, unless otherwise
stated, all our categories are assumed to be $k$-linear. 

A category $\cC$ is called a DG category if for all $X,Y\in \textup{Obj}(\cC)$
$\Hom(X,Y)$ has the structure of a complex of $k$-linear spaces (in other
words, a DG vector space) and the composition maps are associative
$k$-linear maps of DG vector spaces. In
particular, $\Hom(X,X)$ is a DG algebra with a unit.

\begin{ex} \label{DGex}
Given any $k$-linear category $\mathcal{M}$ it is possible to construct a DG category $\cC_{dg}(\mathcal{M})$ with complexes $(M^\bullet,d_M)$ over $\mathcal{M}$ as objects and setting $\Hom(M^\bullet,N^\bullet) = \oplus_n \Hom(M^\bullet ,N^\bullet)_n$, where $\Hom(M^\bullet,N^\bullet)_n$ denotes the component of morphisms of degree $n$, {\it i.e.,} $f_n:M^\bullet\map N^\bullet[n]$ and whose differential is the graded commutator 

\beqn
d(f) &=& d_M\circ f_n - (-1)^nf_n\circ d_N.
\eeqn  
\end{ex}

 Let $\DG$ stand for the category of all small DG categories. The morphisms in
 this category are {\it DG functors}, {\it i.e.,} $F:\cC\to \cC'$ such that for all $X,Y\in
 \textup{Obj}(\cC)$ 

\beqn
F(X,Y):\Hom(X,Y)\map \Hom(FX,FY)
\eeqn is a morphism of DG vector spaces compatible with the compositions and
the units.

\bigskip
\noindent
{\bf The tensor structure:} The tensor product of two DG categories $\cC$ and $\mathcal{D}$ can be defined
in the obvious manner, {\it viz.,} the objects of $\cC\otimes\mathcal{D}$ are
written as $X\otimes Y$, $X\in\textup{Obj}(\cC)$, $Y\in\textup{Obj}(\mathcal{D})$ and one sets

\beqn
\Hom_{\cC\otimes\mathcal{D}}(X\otimes Y,X'\otimes Y') =
\Hom_{\cC}(X,X')\otimes\Hom_{\mathcal{D}}(Y,Y')
\eeqn with natural compositions and units. 

The category of DG functors $\Sheaf{H}om(\cC,\mathcal{D})$
between two DG categories $\cC$ and $\mathcal{D}$ with natural transformations as 
morphisms is once again a DG category. With respect to the above-mentioned tensor product $\DG$
becomes a symmetric tensor category with an inner $\Sheaf{H}om$ functor given by

\beqn
\Hom(\mathcal{B}\otimes\cC,\mathcal{D}) =
\Hom(\mathcal{B},\Sheaf{H}om(\cC,\mathcal{D})).
\eeqn However, in the category of noncommutative spaces (to be defined shortly), this notion of the inner Hom functor needs to be modified.

\bigskip
\noindent
{\bf The derived category of a DG category:} The standard reference for the
construction is \cite{Kel}. We recall some basic facts here. Let $\cC$ be a
small DG category. A right DG $\cC$-module is by definition a DG functor $M:\cC^{\op}\to
\cC_{dg}(k)$, where $\cC_{dg}(k)$ denotes the DG category of complexes of
$k$-linear spaces. Note that the composition of morphisms in the opposite category is defined by the Koszul sign rule: the composition of $f$ and $g$ in $\cC^{\op}$ is equal to
the morphism $(-1)^{|f||g|}gf$ in $\cC$. Every object $X$ of $\cC$ defines
canonically what is called a {\it free} right module $X^\wedge:= \Hom(-,X)$. A morphism of DG modules
$f:L\to M$ is by definition a morphism (natural transform) of DG functors such that $fX:LX\to MX$ is a morphism of complexes for all $X\in\textup{Obj}(\cC)$. We call such an $f$ a
quasiisomorphism if $fX$ is a quasiisomorphism for all $X$, {\it i.e.,} $fX$
induces isomorphism on cohomologies.

\begin{defn} \label{derivedOfDG}
The derived category $D(\cC)$ of $\cC$ is defined to be the localisation of the
category of right DG $\cC$-modules with respect to the class of
quasiisomorphisms. 
\end{defn}

\begin{rem}
With the translation induced by the shift of complexes and triangles coming
from short exact sequences of complexes, $D(\cC)$ becomes a triangulated
category. The Yoneda functor $X\mapsto X^\wedge$ induces an embedding of
$\H^0(\cC)\to D(\cC)$. Here $\H^0(\cC)$ stands for the zeroth cohomology
category whose objects are the same as $\cC$ but the morphisms are replaced by the zeroth
cohomology, {\it i.e.,} $\Hom_{\H^0(\cC)}(X,Y)= \H^0\Hom_{\cC}(X,Y)$. It is also called the homotopy DG category as it produces the homotopy category of complexes over any $k$-linear category $\mathcal{M}$ when specialised to $\cC_{dg}(\mathcal{M})$.
\end{rem}

\begin{defn} \label{pretrDG}
The triangulated subcategory of $D\cC$ generated by the free DG
$\cC$-modules $X^\wedge$ under translations in both directions, extensions and
passage to direct factors is called the {\bf perfect} derived category and
denoted by $\per(\cC)$. A DG category $\cC$ is said to be {\bf pretriangulated}\footnote{Our definition of a pretriangulated category is slightly stronger than \cite{KelDG}, in that, in our definition the homotopy category of such a category is automatically idempotent complete.} if the
above-mentioned Yoneda functor induces an equivalence $\H^0(\cC)\to \per(\cC)$.
\end{defn}

\begin{rem}
A {\it pretriangulated} category does not have a triangulated structure. Rather it is a DG category, which is equivalent to the notion of an
{\it enhanced triangulated category} in the sense of Bondal--Kapranov
\cite{BK}, whose homotopy category is Karoubian. 
\end{rem}

\begin{defn}
A DG functor $F:\cC\to \cD$ is called a Morita equivalence if it induces an
equivalence $F^*:D(\cD)\to D(\cC)$.
\end{defn}
 
\bigskip

\subsection{The category of noncommutative spaces}

\noindent
The definition provided below is a culmination of the works of several people
spanning over two decades including Bondal, Drinfeld, Keller, Kontsevich,
Lurie, Orlov, Quillen and To\"en, amongst others. This list of names is very far from complete and it only reflects the authors ignorance of the history behind this development.

\begin{defn}
The category of noncommutative spaces $\NCS$ is the localisation of $\DG$ with
respect to Morita equivalences.
\end{defn}

Thanks to Tabuada we know that $\DG$ has a cofibrantly generated Quillen model category structure, where the weak equivalences are the Morita equivalences and the fibrant
objects are pretriangulated DG categories. It seems that there was a slight inaccuracy in the proof of the above statement that appeared in \cite{Tab2}, which has now been corrected in \cite{TabThesis}. This enables us to conclude that each object of $\NCS$ has a fibrant replacement, which is a pretriangulated DG category. The tensor product of $\DG$ induces one on $\NCS$ after replacing
any object by its cofibrant model since the tensor product by a cofibrant DG
module preserves weak equivalences. The category $\NCS$ can be regarded as a enhancement of the category of all small idempotent complete triangulated categories. 

\begin{rem}
We have deliberately included  correspondences in the category of
noncommutative spaces. These spaces are somewhat {\it motivic} in nature and it
is expected to be a feature of this geometry. We do not want to treat $\NCS$ as a $2$-category.
\end{rem}

\noindent
However, the inner $\innHom$ functor cannot be derived from $\DG$. Thanks to
To\"en \cite{ToeDG} (also {\it cf.} \cite{KelA-inf}) one knows that there does
exist an inner $\innHom$ functor given by 

\begin{equation} \label{DGintHom}
\Sheaf{H}om(\cC,\cD) = \text{cat. of $A_\infty$-functors $\cC\to\cD$}
\end{equation} here $\cD$ needs to be a pretriangulated DG category which is no
restriction since we know that in $\NCS$ every object has a canonical pretriangulated replacement. The DG structure of $\cD$ endows $\Sheaf{H}om(\cC,\cD)$
with a DG structure as well. We will not be able to discuss
$A_\infty$-categories and $A_\infty$-functors here. Let us mention that a DG category is a special case
of an $A_\infty$-category and we refer the readers to, {\it e.g.,}
\cite{KelA-inf} for a highly readable survey of the same.

\begin{rem} \label{motivicNCS}
The $\Hom$ sets in $\NCS$ are commutative monoids and it is possible to talk about exact sequences in $\NCS$ (see Definition \ref{DGexSeq} below).
\end{rem}

\begin{defn} [Kontsevich]
\end{defn}
{\it \begin{itemize}
\item A noncommutative space (DG category) $\cC$ is called smooth if the bimodule
  given by the DG bifunctor $(X,Y)\mapsto \Hom_{\cC}(X,Y)$ is in
  $\per(\cC^{\op}\otimes \cC)$. 
\item It is called smooth and proper if if is isomorphic in $\NCS$ to a DG algebra whose
  homology is of finite total dimension.
\end{itemize}}

\noindent

\begin{rem}
There is a notion of affinity in this context which just says that a variety is $D$-affine (or derived affine, {\it e.g.,} \cite{BMR} for an analogous notion in the setting of $\D$-modules) if its triangulated category of quasicoherent sheaves is equivalent to the derived category of modules over some (possibly DG) algebra. A theorem of Bondal--Van den Bergh \cite{BonVdb} (see also \cite{SchShi}) asserts that if  $X$ is a quasicompact and quasiseparated scheme, then $D_{Qcoh}(X)$ is equivalent to $D(\Lambda)$ for a suitable DG algebra $\Lambda$ with bounded cohomology. Note that in this theorem $D_{\QCoh}(X)$ denotes the honest derived category of complexes of $\sheaf_X$-modules with quasicoherent cohomologies
and $D(\Lambda)$ likewise. As a consequence we deduce that in the DG setting every proper variety is $D$-affine. 
\end{rem}

\bigskip

\noindent
{\bf Viewing classical geometry in this setting:} We define the {\it DG category of quasicoherent
  sheaves} on an honest scheme $X$ as 

\begin{equation*} \label{ClassGeomIncl}
\cC_{dg} (X) := \cC_{dg}(\QCoh(X)) = \text{DG category of fibrant unbounded complexes over $\QCoh(X)$,}
\end{equation*} which is how we view classical schemes in this framework. It is also
known that $\H^0\cC_{dg} (X) \overset{\sim}{\to} D_{\QCoh} (X)$. As mentioned
above there are reconstruction Theorems available from $\QCoh(X)$ (without any further assumption \cite{Gab,Ros2}) and from $D_{\QCoh} (X)$ (only if the canonical or the anticanonical bundle is ample \cite{BO1} or with the knowledge of the tensor and the triangulated structure \cite{Bal}). They glaringly exclude abelian varieties or (weak) Calabi--Yau varieties, however, for abelian varieties we do have an understanding of the derived category and its autoequivalences \cite{Orl}. 

\begin{rem}
Those who prefer regular holonomic $\D$-modules as substitutes for functions
can perform the above operation after replacing $\QCoh(X)$ by the category of regular holonomic $\D$-modules. 
\end{rem}

Since we have enhanced the morphisms between our spaces by incorporating certain right perfect correspondences, we have also increased the chance of objects becoming
isomorphic. Due to Mukai \cite{Muk} we know that an abelian variety
is derived equivalent to its dual precisely via a {\it correspondence-like}
morphism, which is a Fourier--Mukai transform. Roughly, given any two smooth
projective varieties $X$ and $Y$ and an object in $\Sheaf{E}\in D^b(X\times
Y)$ one constructs an exact Fourier--Mukai transform (also sometimes called an
integral transform) $\Phi^{\Sheaf{E}}_{X\to Y}: D^b(X)\map D^b(Y)$ as follows:

\beqn
\Phi_{X\to Y}^{\Sheaf{E}} (-) = \pi_{Y*}\left( \Sheaf{E}\otimes\pi^*_X
(-)\right),
\eeqn where $\pi_X$ (resp. $\phi_Y$) denotes the projection $X\times Y\to X$
(resp. $X\times Y\to Y$). Here all functors are assumed to be appropriately
derived. The object $\Sheaf{E}$ is called the {\it kernel} of the
Fourier--Mukai transform. In the case of the equivalence between an abelian
variety $A$ and its dual $\hat{A}$ the kernel is given by the Poincar\'e
sheaf $\mathcal{P}$. Given a divisorial correspondence in
$X\times Y$ one can consider the corresponding line bundle on $X\times Y$ and
use that as the kernel of a Fourier--Mukai transform. Conversely, given a
kernel $\Sheaf{E}\in D^b(X\times Y)$ of a Fourier--Mukai transform one obtains
a cycle (correspondence modulo an equivalence relation) in $X\times Y$ by
applying the Chern character to $\Sheaf{E}$. 

\medskip 

\subsection{DG categories up to quasiequivalences} \label{perDG}

We gave a direct method of constructing the category $\NCS$. There is an intermediate notion which one might also want to consider. We call a DG functor $F:\mathcal{C}\map\mathcal{D}$ a {\it quasiequivalence} if the induced maps $\Hom_{\mathcal{C}}(X,Y)\map\Hom_{\mathcal{D}}(FX,FY)$ are quasiisomorphisms for all $X,Y\in\Obj(\mathcal{C})$ and the induced functor $\H^0(F):\H^0(\mathcal{C}\map\H^0(\mathcal{D})$ is an equivalence. The category $\DG$ admits a cofibrantly generated Quillen model category structure whose weak equivalences are quasiisomorphisms \cite{Tab}. Let us denote the homotopy category with respect to this model structure $\Hqe$. Being quasiequivalent is stronger than being Morita equivalent. Therefore, the category of DG categories up to quasiequivalence is bigger (has more non isomorphic objects) than $\NCS$. There is a canonical localisation functor $\Hqe\map\NCS$ inverting the Morita equivalences which are not quasiequivalences, which admits a section functor $\mathcal{A}\map \per_{dg}(\mathcal{A})$, {\it i.e.,} a right adjoint to the canonical localisation functor.

Let us explain the construction of $\per_{dg}(\mathcal{A})$ briefly. For a DG category $\mathcal{A}$ a right $\mathcal{A}$-module, {\it i.e.,} a DG functor from
$\mathcal{A}^{\op}$ to the DG category of complexes over $k$ is called {\it
  semifree} if it admits a countable filtration such that the subquotients are free DG modules (up to shifts), {\it i.e.,} modules formed by arbitrary sums of copies of $\Hom(-,X)$ for some $X\in\Obj(\mathcal{A})$, possibly with shifts. Let us denote the category of
semifree modules over $\mathcal{A}$ by $\textup{SF}(\mathcal{A})$. The
inclusion functor $\textup{SF}(\mathcal{A})\to \mathcal{A}^{\op}$-modules
induces an equivalence of triangulated categories between
$\H^0(\textup{SF}(\mathcal{A}))$ and the derived category of
$\mathcal{A}$ \cite{Dri}. The category $\per_{dg}(\mathcal{A})$ is defined
as the full DG subcategory of $\textup{SF}(\mathcal{A})$ consisting of objects which
become isomorphic to an object in $\per(\mathcal{A})$ after passing on to the
zeroth cohomology category. Roughly speaking, $\per_{dg}(\mathcal{A})$ is
a DG version of $\per(\mathcal{A})$, {\it i.e.,} $\H^0(\per_{dg}(\mathcal{A}))=\per(\mathcal{A})$. 

In fact, the category $\NCS$ is equivalent to the full subcategory of $\Hqe$ consisting of the pretriangulated (or Morita fibrant) DG categories.   

\section{On noncommutative motives}

We begin by reviewing the classical notion of pure motives corresponding to smooth and projective varieties. 

\subsection{Pure motives at a glance} \label{PureMot}

The main steps involved in the construction of effective pure motives from $\Var$ are linearisation, pseudo-abelianisation and finally inversion of the Lefschetz motive, extending the tensor structure of $\Var$ given by the fibre product over $k$. Letting $\sim$ stand for any adequate relation, {\it e.g.,} rational, algebraic, homological or numerical, we define $A^i(X)$ to be the abelian group of algebraic cycles of codimension $i$ in $X$ modulo $\sim$. We define an additive tensor category of correspondences, denoted by $\Corr_\sim$, keeping as objects those of $\Var$ and setting

\beqn
\Corr_\sim(X,Y) &=& \underset{j}{\oplus} A^{\dim\, j}(X\times Y_j),
\eeqn where  each $Y_j$ is an irreducible component of $Y$.

\begin{defn}
An additive category $\mathcal{D}$ is called pseudo-abelian if for any
projector (idempotent) $p\in\Hom(X,X)$, $X\in \Obj(\mathcal{D})$ there exists
a kernel $\ker\, p$.
\end{defn}

There is a canonical {\it pseudo-abelian completion} $\overline{\mathcal{D}}$ of
any additive category $\mathcal{D}$. The objects of $\overline{\mathcal{D}}$
are pairs $(X,p)$, where $X\in \Obj(\mathcal{D})$ and $p\in
\Hom_{\mathcal{D}}(X,X)$ is an arbitrary projector. Define $\Hom$ sets as

\beqn
\Hom_{\overline{\mathcal{D}}}((X,p),(Y,q)) = \frac{\{f\in\Hom_{\mathcal{D}}(X,Y)
\text{\,\, such that\,\,} fp = qf\}}{\{\text{subgroup of\,\,} f \text{\,\, such
  that\,\,} fp=qf=0\}}
\eeqn

We can apply this machinery to construct the pseudo-abelianisation of $\Corr_\sim$. In the resulting category the motive of $\mathbb{P}^n$ decomposes as $\mathbb{P}^n = \textup{pt}\oplus \LL\oplus \LL^{\otimes 2}\oplus\cdots\oplus
\LL^{\otimes n}$. The object $\LL$ is called the {\it Lefschetz motive} and it
should be formally inverted in order to obtain the category of pure
motives and morphisms should also be defined appropriately, but we gloss over these details here. 

Restricting oneself to the subcategory of $\Var$ consisting of connected curves and applying the above-mentioned three steps one obtains the category of motives of curves. This category admits a better description when $\sim$ is chosen to be the rational equivalence relation and morphisms are tensored with $\QQ$.

\begin{prop}[\cite{Man65}] \label{MotCurves}
The category of motives of curves is equivalent to the category of abelian
varieties up to isogeny. 
\end{prop}

\begin{rem}
The functor associates to a curve its Jacobian variety. It turns out that
the category of abelian varieties up to isogeny is abelian and semisimple. 
\end{rem}

The category of motives is expected to be semisimple and Tannakian (Jannsen
  showed that the category of motives modulo numerical equivalence is
  semisimple \cite{Jan}). The category $\NCS$ has some {\it motivic features}: it also has a tensor structure and an inner $\innHom$ functor. But not all objects $T$ are {\it rigid}, {\it i.e.,} the canonical morphism $T\otimes T^\vee\to \Hom(T,T)$ is not an isomorphism for all $T\in\NCS$. However, the
  smooth and proper noncommutative spaces are rigid in the above sense.  

\subsection{Towards noncommutative motives}

The first step of the construction of pure motives entails a linearisation of the
category $\Var$ by including  correspondences. We have argued that
correspondences induce DG functors (indeed, the kernel of a Fourier--Mukai
transform should be thought of as a correspondence). The following Theorem \cite{ToeDG} says that all DG-functors are described by a Fourier--Mukai {\it kernel}, and hence, more relevant to geometry than arbitrary exact functors between triangulated categories.

\begin{thm}[To{\"e}n] \label{DGkernel}
Let $k$ be any commutative ring and let $X$ and $Y$ be quasicompact and separated schemes over $k$ such that $X$
is flat over $\Spec\, k$. Then there is a canonical isomorphism in $\NCS$

\beqn
\cC_{dg}(X\times_k Y) \overset{\sim}{\map} \Sheaf{H}om_c(\cC_{dg}(X),\cC_{dg}(Y)),
\eeqn where $\Sheaf{H}om_c$ denotes the full subcategory of $\Sheaf{H}om$
formed by coproduct preserving quasifunctors, {\it i.e.,} functors between the corresponding zeroth cohomology categories. Moreover, if $X$ and $Y$ are
smooth and projective over $\Spec\, k$, we have a canonical isomorphism in $\NCS$ 

\beqn
\Perf_{dg}(X\times_k Y) \overset{\sim}{\map} \Sheaf{H}om(\Perf_{dg}(X),\Perf_{dg}(Y)),
\eeqn where $\Perf_{dg}$ denotes the full subcategory of $\cC_{dg}$, whose
objects are perfect complexes.
\end{thm}

The above Theorem admits a natural generalisation to abstract DG categories
(not necessarily of the form $\cC_{dg}(X)$ for some scheme $X$), which can also
be found in {\it ibid.}. The above theorem asserts an equivalence of
categories. It can be suitably {\it decategorified}, in order to have an
understanding of the morphisms on the right hand side. 

For DG categories $\mathcal{A},\mathcal{B}$, let us define $\rep(\mathcal{A},\mathcal{B})$ to be the full subcategory of the 
derived category $D(\mathcal{A}^{\op}\otimes\mathcal{B})$ of $\mathcal{A}-\mathcal{B}$-bimodules formed by those $M$, which (under $-\otimes_\mathcal{A} M:D(\mathcal{A})\to D(\mathcal{B})$) send a representable $\mathcal{A}$-module to a $\mathcal{B}$-module, which, in $D(\mathcal{B})$, is isomorphic to a representable $\mathcal{B}$-module. The decategorified statement is that $\Hom(\mathcal{A},\mathcal{B})$ in $\NCS$ is canonically in bijection with the isomorphism classes of objects in
$\rep(\mathcal{A},\mathcal{B})$ {\it ibid.}. If $\mathcal{B}$ is pretriangulated, the objects of
$\rep(\mathcal{A},\mathcal{B})$ are called {\it quasifunctors} as they induce
honest functors $\H^0(\mathcal{A})\to\H^0(\mathcal{B})$. 

Generalising this intuition we conclude that the morphisms in $\NCS$ already
contain all  correspondences. However, $\NCS$ is not an additive category as there is no abelian group structure on the set of morphisms. However, there is a semiadditive structure on $\Hom(\mathcal{A},\mathcal{B})$ given by the direct sum of the kernels of two DG functors or objects in $\rep(\mathcal{A},\mathcal{B})$. We linearise them by passing on to the 
$K_0$-groups of the inner $\innHom$ objects (see, for instance, \cite{KonMot},\cite{Tab2}). 

It is also possible to talk about exact sequences in $\NCS$. We provide one formulation of an exact sequence of DG
categories (see, {\it e.g.,} Theorem 4.11 of \cite{KelDG} for other equivalent
definitions). 

\begin{defn} \label{DGexSeq}
A sequence of DG categories 

\beqn
\mathcal{A}\overset{P}{\map}\mathcal{B}\overset{I}{\map}\mathcal{C}
\eeqn such that $IP=0$ is called exact if and only if $P$ induces an
equivalence of $\per(\mathcal{A})$ onto a thick subcategory of
$\per(\mathcal{B})$ and $I$ induces an equivalence between the idempotent
closure of the Verdier quotient $\per(\mathcal{B})/\per(\mathcal{A})$ and $\per(\mathcal{C})$. 
\end{defn}

\begin{rem}
In the classical setting, if $X$ is a quasicompact quasiseparated scheme, $U\subset X$ a quasicompact open
subscheme and $Z=X\setminus U$, then the following sequence

\beqn
\Perf_{dg}(X)_Z\map \Perf_{dg}(X)\map \Perf_{dg}(U)
\eeqn is exact according to the definition, where $\Perf_{dg}(X)_Z$ denotes the full subcategory of $\Perf_{dg}(X)$
of perfect complexes supported on $Z$. 
\end{rem}

One knows that there is a well-defined $K$-theory functor on $\NCS$, which agrees with Quillen's $K$-theory of an exact
category $\mathcal{B}$, when applied to the Drinfeld quotient of $\cC^b_{dg}(\mathcal{B})$ by its subcategory of acyclic complexes. Now we define a noncommutative analogue of the category of correspondences (a na\"ive version). A more sophisticated approach should treat the category enriched over {\it spectra}, a construction of which can be found in \cite{Tab3}.

\begin{defn} \label{NCmotives}
The category of noncommutative correspondences $\NCCorr$ is the category defined as:

\begin{itemize}
\item $\Obj(\NCCorr) = \Obj(\NCS)$
\item $\Hom_{\NCCorr}(\mathcal{C},\mathcal{D}) = K_0(\rep(\mathcal{C},\mathcal{D}))$
\end{itemize}

\end{defn} 

As a motivation we mention two Theorems: the first Theorem ensures {\it
  linearisation} of $\NCS$, while the second one shows compatibility with localisation.  

\begin{thm}\cite{Dri, DugShi,Tab2} \label{TabuadaSplit}
A functor $F$ from $\NCS$ to an additive category factors through $\NCCorr$ if
and only if for every {\it exact} DG category $\mathcal{B}$ endowed with two
full exact DG subcategories $\mathcal{A},\mathcal{C}$ which give rise to a
{\it semiorthogonal decomposition}
$\H^0(\mathcal{B})=(\H^0(\mathcal{A}),\H^0(\mathcal{C}))$ in the sense of
\cite{BO2}, the inclusions induce an isomorphism $F(\mathcal{A})\oplus
F(\mathcal{C})\overset{\sim}{\map} F(\mathcal{B})$. 
\end{thm}

\noindent
Such a functor is called an {\it additive invariant} of noncommutative
spaces. The simplest example is $\mathcal{A}\longmapsto
K_0(\per(\mathcal{A}))$.

\begin{thm} \cite{DugShi}
The functor $\mathcal{A}\longmapsto K(\mathcal{A})$ (Waldhausen $K$-theory) yields, for each short
exact sequence $\mathcal{A}\to \mathcal{B}\to\mathcal{C}$ in $\NCS$, a long
exact sequence
\beqn
\cdots\map K_i(\mathcal{A})\map K_i(\mathcal{B})\map
K_i(\mathcal{C})\map\cdots \map K_0(\mathcal{B})\map K_0(\mathcal{C}).
\eeqn

\end{thm}

\begin{rem}
The category $\NCCorr$ is additive and the composition is induced by that of
$\NCS$. Certain non-isomorphic objects of $\NCS$ become isomorphic in $\NCCorr$, {\it
  e.g.,} it is shown in \cite{KelLocDG} that each finite dimensional algebra of finite global dimension becomes isomorphic to a product of copies of $k$ in $\NCCorr$, whereas such a thing is
true in $\NCS$ if and only if the algebra is semisimple.
\end{rem}

 We should perform a formal idempotent completion (or pseudo-abelian completion) of $\NCCorr$ as discussed in Subsection
\ref{PureMot} in order to obtain the category of noncommutative motives, which is denoted by $\NCM$. It follows from Be{\u\i}linson's description of the derived category of coherent sheaves on $\mathbb{P}^n$ \cite{Bei} and the above remark that $\cC_{dg}(\mathbb{P}^1) \simeq \cC_{dg}(\mathbb{A}^1)\oplus\cC_{dg}(\pt)$ is also isomorphic in $\NCCorr$ to $\cC_{dg}(\pt)\oplus\cC_{dg}(\pt)$, whence $\cC_{dg}(\mathbb{A}^1)\simeq \cC_{dg}(\pt)$, {\it i.e.,} the Lefschetz motive is isomorphic to the identity element.

\medskip

A careful reader should have noticed that we have glossed over the issue of the
choice of the equivalence relation, which was central to the
construction of the category of pure motives in the classical setting. Manin mentioned in \cite{Man65}
(end of Section 3) that every  cohomology theory should be a
cohomological functor on the category of $\Corr_\sim$, {\it i.e.,}
every correspondence in $\Corr_\sim(X,Y)$ should induce a well-defined morphism
$H^*(X)\to H^*(Y)$. Now we turn the argument around. Elements of $\rep(\mathcal{C},\mathcal{D})$ induce morphisms between $\mathcal{C}$ and $\mathcal{D}$. Our spaces are defined in terms of the (quasicoherent) cohomologies that they admit. Mostly cohomology theories appear as cohomology groups of a certain canonically defined cochain complex satisfying a bunch of axioms. We pretend that a morphism (a functor) in $\NCM$ is a morphism between
the cohomology theories on the two spaces, as if given by some correspondence. If the question about universal cohomology theory is resolved, then probably one would like to argue that the elements of $\rep(\mathcal{C},\mathcal{D})$ are the ones which induce distinct morphisms between their universal cohomologies. If that turns out to be false then one can call an equivalence relation {\it universal} if it identifies two
correspondences which induce isomorphic morphisms between the corresponding
{\it universal cohomology theories} and then consider correspondences modulo this equivalence relation. Note that in $\NCM$ we set the
Grothendieck group of $\rep(\mathcal{A},\mathcal{B})$ as morphisms between
$\mathcal{A}$ and $\mathcal{B}$. Chow correspondences are obtained by taking the rational equivalence relation. The connection should be an analogue of the {\it Chern character} map which identifies the K-theory with the Chow group after tensoring with $\QQ$. The readers are referred to \cite{BGNT} for a possibly relevant treatment of the Chern character. 

\subsection{Motivic measures and motivic zeta functions} \label{MotMeasure}

We present a rather simplistic point of view on motivic measures. With respect
to a motivic measure it is possible to develop a theory of {\it motivic
  integration} (see, {\it e.g.,} \cite{DL}), which we shall not discuss
here. This technology was invented by Kontsevich drawing inspiration from the works of Batyrev. A useful and instructive reference is, {\it e.g.,} \cite{Loo}. 

Let $\Sch_k$ be the category of reduced schemes of finite type (or reduced varieties) over $k$. Consider the Grothendieck ring of $\Sch_k$, denoted by
$K_0(\Sch_k)$, which is defined as the free abelian group generated by isomorphism classes of objects
in $\Sch_k$ modulo relations (often called scissor-congruence relations)

\begin{equation} \label{excisionCond}
[X] = [Z] + [X\setminus Z],
\end{equation} where $Z$ is a closed subscheme of $X$. The multiplication is given by
the fibre product over $k$. There is a unit given by the class of $\Spec\, k$. 

Every $k$-variety admits a finite stratification $X=X^0\supset X^1\supset\cdots\supset X^{d+1} = \emptyset$ such that $X^k \setminus X^{k+1}$ is smooth. Moreover, any two such stratifications admit a common refinement. Therefore $[X]=\sum_k [X^k\setminus X^{k+1}]$ is unambiguously defined and, in fact, it can be shown that $K_0(\Sch_k)$ is generated by complete and nonsingular varieties. The structure of $K_0(\Sch_k)$ as a ring is rather complicated. It is known that it is not an integral domain \cite{Poo}. However, it admits interesting ring homomorphisms to some rings, which turn out to be quite useful in various cases.

Let $A$ be any commutative ring. An $A$-valued {\it motivic measure} is a ring
homomorphism $\mu_A:K_0(\Sch_k)\to A$. We write $\mu = \mu_A$ if there is no chance of confusion. If $A$ has a unit the homomorphism is required to be unital.

\begin{ex}
Let $k=\CC$, $A=\ZZ$ and $\mu(X) = \chi_c(X)$, {\it i.e.,} the Euler characteristic with compact supports.
\end{ex}

\begin{ex}
Let $k=\CC$, $A=K_0(\textup{HS})$, {\it i.e.,} the Grothendieck ring of Hodge structures and $\mu(X) = \chi_h(X)$ such that 
\beqn
\chi_h(X) = \sum_r (-1)^r[\H_c^r(X,\QQ)] \in K_0(\textup{HS}),
\eeqn which is called the Hodge characteristic of $X$.
\end{ex}

\begin{ex} \label{HasseWeilmeasure}
Let $k=\mathbb{F}_q$, $A=\ZZ$ and $\mu(X) = \# X(\mathbb{F}_q)$, {\it
  i.e.,} the number of $\mathbb{F}_q$-points.
\end{ex}

Let us fix an $A$-valued motivic measure $\mu$ and, for a smooth $X\in\Sch_k$,
let $X^{(n)}$ denote the $n$-fold symmetric product of $X$. Set $X^{(0)}:=
\Spec\, k$. Then associated to
$\mu$ there is a {\it motivic zeta function} (possibly due to Kapranov
\cite{KapZeta}) of $X$ defined by the formal series

\begin{equation}
\zeta_\mu(X,t) = \sum_{n=0}^{\infty}\mu\left(X^{(n)}\right)t^n \in A[[t]].
\end{equation}

\begin{ex}
If $k=\mathbb{F}_q$, $A=\ZZ$ and $\mu(X) = \# X(\mathbb{F}_q)$ as in Example
\ref{HasseWeilmeasure} one recovers the usual Hasse--Weil zeta function of
$X$. Indeed, the $\mathbb{F}_q$-valued points of $X^{(n)}$ correspond to the
effective divisors of degree $n$ in $X$. 
\end{ex}

Let us denote $\mu(\Aff_k^1)$ by $\LL$. Then we have the following rationality statement for curves
(see Theorem 1.1.9 {\it ibid.}).

\begin{thm}
If $X$ is any one dimensional variety (not necessarily non-singular) of genus
$g$, then $\zeta_\mu(X,t)$ is rational. Furthermore, the rational function
$\zeta_\mu(X,t)(1-t)(1-\LL t)$ is actually a polynomial of degree $\leqslant
2g$ and satisfies the functional equation below.

\begin{equation}
\zeta_\mu(X,1/\LL t) = \LL^{1-g}t^{2-2g}\zeta_\mu(X,t)
\end{equation}

\end{thm}

\begin{rem}
The rationality statement fails to be true in higher dimensions, {\it e.g.,}
if $X$ is a complex projective non-singular surface of  genus
$\geqslant 2$ \cite{LarLun1}. In fact, a complex surface $X$ has a rational
motivic zeta function if and only if it has {\it Kodaira dimension} $-\infty$
\cite{LarLun2}.
\end{rem} 

\subsection{Noncommutative Calabi--Yau spaces}

This section attempts to introduce zeta functions of noncommutative curves
{\it in a motivic framework} and possibly extract some arithmetic information out of
them. That the zeta functions of varieties contain crucial arithmetic information is a gospel truth by now. 

Before we move forward let us mention that such ideas are prevalent in
noncommutative geometry, {\it e.g.,} Connes' spectral realisation of the zeros
of the Riemann zeta function \cite{ConRiemann1,ConRiemann2}. Some other
important works in this direction are
\cite{ConMarRam},\cite{Den1,Den2},\cite{Pla} and \cite{HaPau}, to mention only a few. Also the
readers should take a look at \cite{Mar} for a more holistic point of view.

\medskip

Following Proposition \ref{MotCurves} we argue that the category of noncommutative
motives of noncommutative curves should be equivalent to the full subcategory of $\NCM$
generated by DG categories which resemble those of abelian varieties, {\it
  i.e.,} the inclusion of abelian varieties inside $\NCM$ (see Equation
\eqref{ClassGeomIncl}). Given an abelian surface the cokernel of the
multiplication by $2$ map (isogeny) is a Kummer surface with $16$ singular
points, whose (minimal) resolution of singularities is a $K3$ surface. It is
an example of a Calabi--Yau manifold of dimension $2$. So even if we look at
motives of curves Calabi--Yau varieties show up rather naturally. We propose
to treat such varieties as they are, rather than working up to isogenies. Calabi--Yau varieties are interesting from the point of view of physics as well. We assume that a
Calabi--Yau variety is just a variety, whose canonical class is trivial (no
assumption on the fundamental group). 

In a $k$-linear category $\mathcal{A}$ an additive
  autoequivalence $S$ is called a Serre functor if there exists a bifunctorial
  isomorphism $Hom(A,B)\overset{\sim}{\to} Hom(B,SA)^*$ for any two $A,B\in\Obj(\mathcal{A})$. If
it exists it is unique up to isomorphism. If $X$ is a smooth projective variety of dimension $n$,
the Serre functor is given by $(-\otimes\omega_X)[n]$, where $\omega_X$ is the canonical sheaf of $X$. 
The existence of a Serre functor corresponds to that of Grothendieck--Serre duality.

\begin{defn}
A DG category $\mathcal{C}$ in $\NCS$ is called a noncommutative Calabi--Yau
space of dimension $n$ if $\H^0(\mathcal{C})$ is triangulated ({\it i.e.,}
$\mathcal{C}$ is pretriangulated as in Definition \ref{pretrDG}) with the finiteness condition $\sum_p \dim_k \,\Hom_{\H^0(\mathcal{C}}(X,Y[p]) <\infty$ for all $X,Y\in\Obj(\H^0(\mathcal{C})$, and if there exists a natural isomorphism between the Serre functor and $[n]$. In other
words, there exists bifunctorial isomorphisms
$\Hom(A,B)\overset{\sim}{\to}\Hom(B,A[n])^*$ in $\H^0(\mathcal{C})$. 
\end{defn} 

Kontsevich originally defined a noncommutative Calabi--Yau space as a small triangulated category satisfying the strong finiteness condition mentioned above, with an isomorphism between the Serre functor and $[n]$. We have enhanced it to the DG level. It should be bourne in mind that the homotopy category of a pretriangulated category is  idempotent complete. It is clear that if $X$ is a Calabi--Yau variety then $\cC_{dg}(X)$ is a Calabi--Yau space in the above sense. Purely at the triangulated level there are other interesting examples of Calabi--Yau spaces of dimension $2$ arising from quiver representations and commutative algebra {\it cf.} Section 4 of \cite{KelRei}. When such a triangulated Calabi--Yau category of dimension $d$ is endowed with a {\it cluster tilting subcategory} it is possible to construct a Calabi--Yau DG category (in the above sense) of dimension $d+1$ \cite{TabCY}. 

Let us denote by $\MC$ the full additive subcategory of $\NCM$ consisting of noncommutative Calabi--Yau spaces. 

\begin{ex}
It is expected that via a noncommutative version of the construction of the
Jacobian of a curve the category of motives of noncommutative curves can be
seen as a full subcategory of $\MC$. The way to view an abelian variety in this
setting is not clear to the author yet. The category $\MC$ contains honest
elliptic curves (as they are their own Jacobians) as given by the inclusion of
classical geometry in this setting (see Equation \eqref{ClassGeomIncl}). The
noncommutative torus $\TT_\theta^\tau$ is also included via its DG derived
category of holomorphic bundles. It is isomorphic to $\cC_{dg}(X_\tau)$, where $X_\tau = \CC/(\ZZ +\tau\ZZ)$, via a Fourier--Mukai type functor (see Proposition 3.1 \cite{PolSch}). J. Block \cite{Blo1,Blo2} suggests a more conceptual framework for such dualities to exist. The rough idea is to construct a differential graded algebra from a complex torus $X$ associated to a deformation parameter in $\H\H^2(X)$ and look at the DG category $DG(X)$ of {\it twisted} complexes, {\it i.e.,} DG modules over that algebra equipped with a super connection compatible with the differential of the algebra. One can construct a (curved) differential graded algebra corresponding to the dual torus as well with a curvature contribution given by the deformation parameter, whose DG category of twisted complexes will be quasiequivalent to $DG(X)$ via some sort of a deformed Poincar\'e bundle (essentially a correspondence).  
\end{ex}

\medskip

\subsection{The motivic ring of $\NCS$}\footnote{The author was kindly notified by B. Keller and the referee that the motivic ring constructed in \cite{BLL}, which is what we describe here, is actually isomorphic to the $0$ ring \cite{Tab2}. However, with some appropriate finiteness conditions thrown in, this problem has been rectified, ({\it e.g.,} \cite{TabThesis}).}

Let us recall from Section \ref{MotMeasure} that an $A$-valued motivic measure
$\mu$ is a ring homomorphism from $K_0(\Sch_k)\to A$. We have replaced the category
of $k$-schemes by the category of noncommutative spaces $\NCS$. We need an appropriate notion of the Grothendieck ring of $\NCS$, which we would like to call the motivic ring of $\NCS$. 

Since every object in $\NCS$ is quasiequivalent to a pretriangulated
DG category we seek a Grothendieck ring of pretriangulated DG categories. In \cite{BLL} the authors precisely construct a Grothendieck ring of
pretriangulated DG categories, which is essentially the Grothendieck ring of $\Hqe$. It was pointed out by the authors that it is
crucial to work with DG categories (and not honest triangulated ones) as the tensor
product of two triangulated categories does not have a natural triangulated structure in general. Let us
briefly recall their construction. 

The Grothendieck ring $\mathcal{G}$ is generated as a free abelian group by
the isomorphism classes of pretriangulated DG categories in $\NCS$ (or
quasiequivalence classes of objects in $\DG$) modulo
relations analogous to those of $K_0(\Sch_k)$. The authors reinterpret the excision 
relations as those coming from {\it semiorthogonal decompositions} (see \cite{BO2} for the details of
semiorthogonal decomposition). One writes $[\mathcal{B}] = [\mathcal{A}] +
[\mathcal{C}]$ if and only if there exist representatives $\mathcal{A}'$,
$\mathcal{B}'$, $\mathcal{C}'$ in $[\mathcal{A}]$, $[\mathcal{B}]$,
$[\mathcal{C}]$ respectively such that 

\begin{enumerate}
\item $\mathcal{A}'$, $\mathcal{C}'$ are DG subcategories of $\mathcal{B}'$,
\item $\H^0(\mathcal{A}')$, $\H^0(\mathcal{C}')$ are admissible subcategories
  of $\H^0(\mathcal{B}')$,
\item \label{semiDecomp} $(\H^0(\mathcal{A}'),\H^0(\mathcal{C}'))$ is a semiorthogonal
  decomposition of $\H^0(\mathcal{B}')$. 
\end{enumerate}

\begin{rem}
Part \eqref{semiDecomp} implies that
$\H^0(\mathcal{A}')=(\H^0(\mathcal{C}'))^\perp$, which is Lemma 2.25 in
\cite{BLL}. An exact sequence $\mathcal{A}\map\mathcal{B}\map\mathcal{C}$ of pretriangulated DG categories ({\it cf.}
Definition \ref{DGexSeq}) induces an {\it exact} sequence of honest
triangulated categories
$\H^0(\mathcal{A})\map\H^0(\mathcal{B})\map\H^0(\mathcal{C})$ by
definition. However, existence of a semiorthogonal decomposition is a stronger
condition. It says that $\H^0(\mathcal{C})$ is a triangulated subcategory of
$\H^0(\mathcal{B})$ and $\H^0(\mathcal{A})=(\H^0(\mathcal{C}))^\perp$, {\it
  i.e.,} the sequence is {\it split} ({\it cf.} Theorem \ref{TabuadaSplit}). It is plausible that one obtains
something sensible by allowing all possible exact sequences as relations. 
\end{rem}

\noindent
The product $\bullet$ is defined as follows:

\beqn
\mathcal{A}_1\bullet \mathcal{A}_2 := \per_{dg}(\mathcal{A}_1\otimes
\mathcal{A}_2),
\eeqn where $\per_{dg}(\mathcal{A})$ is a pretriangulated DG category
described in subsection \ref{perDG}.

\noindent
The product $\bullet$ preserves quasiequivalences of DG categories and hence
descends to a product on $\mathcal{G}$. It is proven in \cite{BLL} that the
product is associative and commutative. There is a unit given by the class of
$\cC^b_{dg}(k)$, {\it i.e.,} the DG category of finite dimensional chain complexes
over $k$. That this product corroborates the fibre product of varieties is justified by Theorem 6.6 {\it ibid.}.

\medskip

\begin{rem}
The motivic ring of $\NCS$ should be the above-mentioned ring with quasiequivalences replaced by Morita equivalence. There will be a canonical ring homomorphisms corresponding to the localisation functor $\Hqe\map\NCS$. Since a quasiequivalence is also a Morita equivalence the ring homomorphism will be surjective identifying elements which are Morita equivalent but not quasiequivalent.
\end{rem}

It follows that the product of two noncommutative Calabi--Yau categories is again a
noncommutative Calabi--Yau category, {\it i.e.,}  $\mathcal{A}\bullet \mathcal{B}$ is a noncommutative
Calabi--Yau DG category of dimension $m+n$ for $\mathcal{A},\mathcal{B}\in\MC$ of dimensions $m$, $n$ respectively. Indeed, the finiteness condition follows from K\"unneth formula. To check the existence of the Serre functor first observe that

\beqn
&&\Hom_{\H^0(\cA\otimes\cB)}(A\otimes B,A'\otimes B'[m+n]) \\
&=& \H^0\Hom_{\cA\otimes\cB}(A,A')\otimes\Hom_{\cA\otimes\cB}(B,B')[m+n]\\
&=& \H^{m+n}_{\cA\otimes\cB}(\Hom(A,A')\otimes\Hom(B,B'))\\
&=& \H^{m+n}\left((\oplus_i
\Hom^i_{\cA\otimes\cB}(A,A')\otimes\Hom^{k-i}_{\cA\otimes\cB}(B,B'))^\bullet\right) \\
&=& \oplus_l \Hom_{\H^0(\cA\otimes\cB)}^{l+m}(A,A')\otimes
\Hom_{\H^0(\cA\otimes\cB)}^{n-l}(B,B') \\
&=& \oplus_l
\Hom^l_{\H^0(\cA\otimes\cB)}(A,A'[m])\otimes\Hom^{-l}_{\H^0(\cA\otimes\cB)}(B,B'[n])\\
&=& \oplus_l
\Hom^l_{\H^0(\cA\otimes\cB)}(A',A)^*\otimes\Hom^{-l}_{\H^0(\cA\otimes\cB)}(B',B)^*\\
&=& \left(\oplus_l \Hom^l_{\H^0(\cA\otimes\cB)}(A',A)\otimes
\Hom^{-l}_{\H^0(\cA\otimes\cB)}(B',B)\right)^*\\
&=& \H^0(\Hom^\bullet_{\cA\otimes\cB}(A'\otimes B',A\otimes B))^*\\
&=& \Hom_{\H^0(\cA\otimes\cB)}\Hom(A'\otimes B',A\otimes B)^*
\eeqn This proves that $\H^0(\cA\otimes\cB)$ has the right Serre functor $[n+m]$. 

Now, it follows from \cite{LyuMan} that the existence of the Serre functor $[n]$ is equivalent to the isomorphism 
$\Hom_k(\cA(-,?),k) \simeq \cA(?,-[n])$ in $D(\cA\otimes \cA^{op})$. Since $\cA \mapsto \per_{dg}(\cA)$ is a monoidal Morita isomorphism this equivalence of bimodules is preserved under the functor.

\vspace{2mm}


A $\G$-valued motivic measure (the ring homomorphism being the identity map) is a universal one. Let us denote the class of $X\in\NCS$ inside $\G$ by $[X]$. Given any variety $Y$ we know that $\mathcal{C}_{dg}(Y\times Y)$ is quasiequivalent (in particular, Morita equivalent) to $\Sheaf{H}om_c(\mathcal{C}_{dg}(Y),\mathcal{C}_{dg}(Y))$. Let us define inductively $X^n = \Sheaf{H}om_c(X^{n-1},X)$ and $X^1 = X$. Then the universal $\G$-valued motivic zeta function of $X\in\NCS$ is given by 

\beqn
\zeta_{\mu_{\G}}(X,t) = 1 + \sum_{n=1}^\infty [X^n]
t^n \qquad \in \G[[t]].
\eeqn 

It is shown in \cite{BLL} that there is a canonical surjective ring
homomorphism $K_0(\Sch_k)\to \G_{hon}$ with $(\LL-1)$ in the kernel, where
$\G_{hon}$ is the subring of $\G$ generated by certain pretriangulated DG categories associated to honest smooth projective varieties over $k$.  

\medskip

\begin{rem}
Since the DG category of holomorphic bundles on the noncommutative torus $\TT_\theta^\tau$ is equivalent to $\cC_{dg}(X_\tau)$, where $X_\tau = \CC/(\ZZ +\tau\ZZ)$ \cite{PolSch}, the zeta functions of $\TT_\theta^\tau$ and $X_\tau$ are the same. However, the equivalence is given by a non-trivial Fourier--Mukai type (correspondence) functor. The B-model of a conformal field
theory associates to a complex torus its derived category of coherent
sheaves. The equivalence perhaps indicates that deforming the complex
torus to a noncommutative torus does not produce anything new for the B-model. 
\end{rem}

One nagging point is that certain natural topological constructions do not
allow us to define a category in which the composition of morphisms obeys
associativity (it is associative only up to homotopy). Hence some mathematicians
have resorted to working with $A_\infty$ categories which encode such
properties, {\it e.g.,} \cite{KonSoiAinf}. The world of $A_\infty$ categories subsumes that of DG
categories. However, one knows that every $A_\infty$ category is
quasiisomorphic to a DG category in a functorial manner ({\it e.g.,}
\cite{CosTCFT}).


\bibliographystyle{abbrv}
\bibliography{/Users/smahanta/Professional/math/MasterBib/bibliography}

\end{document}